\documentclass[a4paper,10pt,DIV12]{article}
\usepackage{indentfirst}
\usepackage{latexsym}
\usepackage{amssymb, amsmath,amsfonts}
\usepackage{manfnt}
\usepackage{mathrsfs}
\usepackage[ansinew]{inputenc}

\usepackage{epsfig, color}

\newcommand{\eproof}{\hspace*{ \fill } $ \Box $ \vspace{ 0.3cm }}
\newenvironment{proof}{ \emph{Proof}. }{ \eproof }

\title{The fine structure of 321 avoiding involutions.}
\author{Piera Manara \thanks{Università di Parma ({\tt piera.manara@fis.unipr.it})} and Claudio Perelli Cippo \thanks{
 Politecnico di Milano
 ({\tt claudio.perelli\_cippo@polimi.it}}.}
\begin{document}
\maketitle
In questo lavoro si studiano proprietà delle involuzioni appartenenti alla classe $Av(321)$. Si calcolano le funzioni generatrici algebriche dell'insieme delle involuzioni di $Av(321)$ e di alcuni suoi sottoinsiemi. Precisamente si calcolano le funzioni generatrici algebriche delle involuzioni espansione di 12, delle involuzioni espansione di 21, delle involuzioni semplici e delle involuzioni espansione di semplici appartenenti alla classe $Av(321)$.\\
Si caratterizzano i grafici delle involuzioni semplici.\\ Infine si dà un{'} interpretazione combinatoria di alcuni dei risultati ottenuti, per mezzo di opportune classi di cammini di Motzkin.\\

\textbf{Abstract}.  We study the involutions belonging to the class of  321 avoiding permutations. We calculate the algebraic generating functions of the set containing the involutions avoiding 321 and of some of its subsets. Precisely we determine the algebraic generating functions of the involutions that are expansions of 12, of those expansions of 21, of the simple ones and of their expansions.\\ The graphics of the simple involutions are caracterized. Being the simple involutions avoiding 321 counted by Riordan's numbers, a combinatoric interpretation of the results is illustrated through a class of Motzkin paths. Another interpretation is given through Dyck paths.

\section{Introduction.}

The aim of this work is to study the set $I(321)$, consisting of all involutions in the class of  321 avoiding permutations{,} and some of its subsets. The results are obtained by means of the \textit{substitution decomposition} properties, particularly of the \textit{involution decomposition},
{given in \cite{a},  \cite{B} and \cite{MR}, and the techniques of generating functions, as used in \cite{a}, \cite{B}, \cite{AA}.}\\
We briefly enumerate the theorems we use and some definitions, maintaining almost always the terminology used in the cited papers, to which we refer for the demonstrations and some well known basic definitions, as the ones of permutation, class of permutations, inflation of a permutation. \\

For a permutation set $S$, we denote by $S_{n}$ the set of the permutations in $S$ of length $n$, and we refer to $f(x)\, =\,\sum |S_{n}| x^{n}$ as the \textit{generating function for S}.\\

An \textit{interval } in the permutation $\pi$ is a set of contiguous indices $\cal I$ = $[a,\,b]$, such that the set of values $\pi({\cal I }) = \left\{\pi(i):\,i \in {\cal I } \right\}$ is also contiguous.\\

A permutation $\pi \in S_{n} $ is said to be \textit{simple } if it contains only the intervals $0,\,1,\,[1,\ldots{,}n]$. \\
The permutations 1,\,12,\,21 are simple; there is no simple permutation for $n=3$.\\

The permutations \textit{avoiding the pattern 321} {constitute the} class, $Av(321)$, such that any $\pi \in Av(321)$
contains no descending sequence{s} of length $n\geq 3$.\\
Interesting properties of the graphic of a 321 avoiding permutation are given in \cite{a}, {Section} 2;
in particular such a graphic is the merge of two increasing sequences.\\

An \textit{involution} is a permutation {$\pi$ such that $\pi(\pi(i))=i$ for all $i=1,\ldots ,n$}.\\
The graphic of any involution has obviously a symmetry with respect to the line $y=x$. So, the graphic of an involution avoiding $(321)$
 shows two symmetric increasing sequences.\\

The properties regarding \textit{substitution decomposition} and \textit{involution decomposition}, which we call \textit{structure theorems}, are given in \cite{a}, \cite{B} e \cite{AA}; we enumerate in the following those propositions because fundamental for our work.\\
\\
{\bf Proposition 1.1} (See Albert and Atkinson \cite{AA}.) \textit{Every permutation, except 1, is the inflation of a unique simple permutation of length at least 2}.\\
{(This means that every permutation $\pi$ determines the unique simple permutation of which $\pi$ is an inflation).}\\
\\
{\bf Proposition 1.2} (See \cite{AA}.) \textit{If $\pi\,=\,\sigma[\alpha_{1},\ldots,\alpha_{m}]$ where $\sigma$ is simple of length $m\geq 4$, then the $\alpha_{i}'s$ are unique.}\\

In the following we shall also say  shortly  \textit{of type} 12 and \textit{of type} 21 for
the substitutions which are inflations respectively of 12 and 21 \textit{(sum decomposable, minus decomposable)}.\\
As for the case of the simple permutations 12 and 21 a unique decomposition can still be { canonically} given, as in the following.\\
\\
{\bf Proposition 1.3} (See \cite{AA}.)\textit{ If $\pi$ is an inflation of 12, then there is a unique $\alpha_{1}$, not of type 12, such that $\pi\,=\,12[\alpha_{1},\alpha_{2}]$  for some $\alpha_{2}$, which is itself unique. The same holds with 12 replaced by 21, and $\alpha_{1}$ not of type 21.}\\

In \cite{B} the structure properties are applied specifically to the involutions, which we study in this work. We recall the following propositions.\\
\\
{\bf Proposition 1.4} \textit{{A permutation } $\pi\,=\,12[\alpha_{1},\,\alpha_{2}],$ is an involution if and only if $\alpha_{1}$ and $\alpha_{2}$ are involutions.}\\
\\
{\bf Proposition 1.5} \textit{Let be $\pi = \sigma[\alpha_{1},\ldots,\alpha_{m}]$ inflation of a simple permutation $\sigma \neq 21$. Then if $\pi$ is an involution, also $\sigma$ is an involution, and the following equalities hold: $\alpha_{i}=\alpha^{-1}_{\sigma^{-1}(i)}=\alpha^{-1}_{\sigma(i)}$, for $i=1,\ldots,m$}.\\
(The equalities say that every transposition in $\sigma$ {must} be inflated with a couple of substitutions one the inverse of the other).\\
\\
{\bf Proposition 1.6}  \textit{The involutions that are inflations of 21 are precisely those of the form:\\
 $21[\alpha_{1},\,\alpha_{2}]$, where neither $\alpha_{1}$ nor $\alpha_{2}$ are of type $21$, and $\alpha_{1} =\alpha^{-1}_{2}$;\\
$321[\alpha_{1},\,\alpha_{2},\, \alpha_{3}]$, where neither $\alpha_{1}$ nor $\alpha_{3}$ are of type $21$,  $\alpha_{1} =\alpha^{-1}_{3}$, and $\alpha_{2}$ is an involution.}\\

 \section{Properties of the involutions in $Av(321)$.}

We denote $I(321)$ the set of the involutions in  $Av(321)$, and $I(321)_{n}$ its subset of the involutions of length $n$.\\
{ Recall that the left to right \textit{maxima }of a permutation are those elements wich dominate all of their predecessors.\\
We give in the following some results about involutions in $I(321)$.}\\
\\
{\bf Proposition 2.1} \textit{If $\pi\in I(321)$, $\pi\,=\,12[\alpha_{1},\,\alpha_{2}],$ then $\alpha_{1}\, {\rm and}\, \alpha_{2}$ are in  I(321).}\\
\\
{\bf Proposition 2.2} \textit{If $\pi\in I(321)$ is of type $21$, then $$\pi\,= \,21[1\,2\,\ldots\, m,\, 1\,2\,\ldots \,m]\,= \,(m+1)\,(m+2) \ldots\, (2m)\, 1\,2 \ldots\, m,$$
and $\pi$ has even length, $\pi \in I(321)_{2m}$.}\\

\begin{proof} Let $\pi = 21[\alpha_{1},\, \alpha_{2}] $ be an involution; by Proposition 1.3, $\alpha_{1} =\alpha^{-1}_{2}$, then $\pi$ has length $n=2m$. Being  $\pi $ of type 21, its maxima
are only in the left block; because $\pi$ is in  Av(321), in the right block one has the ascending sequence $\alpha_{2}\,=\,1\,2 \,\ldots \,m.$ Hence one derives $\alpha_{1}\,= \,\alpha^{-1}_{2} =\,1\,2\,\ldots\, m$, and the thesis.\end{proof}

 Recall that an involution $\pi\in S_{n}$ can be written as a product of cycles of length 1 and 2, respectively {\it fixed points} and {\it transpositions}:
$$\pi\,=\,(m_1,M_1)(m_2,M_2)\ldots(m_h,M_h),\,\,{\rm con}\quad m_i\leq M_i,\, i=1,\ldots,h .$$
An involution $\pi\in S_{2n+1}$ has an odd number greater or equal to $1$ of fixed points; $\pi\in S_{2n}$ has no fixed points or an even number of them.\\
\\
{\bf Proposition 2.3} \textit{A necessary and sufficient condition for  an involution $\pi$ being in I(321) is to hold the inequalities $1=m_1<m_2<\cdots <m_h$ and $M_1<M_2<\cdots <M_h$ { (so that the $M_i$ are left to right maxima of $\pi$)}.}\\
\\
{\bf Proposition 2.4} \textit{If $\sigma\in I(321)_n$, $n>2$, is a simple involution, then $n$ is even. If  $\pi$ is an inflation of a simple involution $\sigma \in I(321)_n$,  $n>2$, then the length of  $\pi$ is even.}\\

 \begin{proof} A simple involution $\sigma \in I(321)$ has no fixed points, because of the symmetry of
the graphic of $\sigma$, hence of the two ascending sequences. But because any involution of odd length must have an odd number of fixed points, the length of $\sigma \in I(321)$ cannot be  odd.\\
As for an inflation of a simple involution, it must have an even length because always obtained
through couples of substitutions one the inverse of the other, by Proposition 1.5.\end{proof}
\\(For example, \textbf{35}1{\bf 6}24[12,1,12,1,1,1] = \textbf{457}12 \textbf{8}36.)\\
\\
{\bf Proposition 2.5} \textit{The set $I(321)$ has  infinitely many simple involutions.}\\

\begin{proof} { Indeed it is immediate to see that, for every $n\geq 3$,
the involution $\sigma \in I(321)_{4n-6}$ defined as
$\sigma \,=\;n\,(n+2)\ldots (n+2(n-2))\,1\,(n+2(n-2)+1)\,2\,\ldots (4n-6)(n-1)((n-1)+2)\ldots((n-1)+2(n-2))\;$
the sequence of whose maxima is $\left\{\,n, \,n+2,\,\ldots,\, n+2(n-2),n+2(n-2)+1,\,\ldots,n+2(n-2)+n-2\,\right\}$,
is a simple involution.\\
(For instance  $\;\textbf{35}1{\bf 6}24\;=\,(1,3)(2,5)(3,6) \in I(321)_{6}$, \\
$\textbf{468}1\textbf{9}2 (\textbf{10})357\,=\,(1,4)(2,6)(3,8)(5,9)(7,(10))\;\in
I(321)_{10},$\\
$\;\textbf{579(11)}1 \textbf{(12)}2 \textbf{(13)}3 \textbf{(14)}468(10)\;\in
I(321)_{14},\,\dots$)}\end{proof}
\\
{\bf Proposition 2.6} \textit{{ Let} $\sigma \,\in \,I(321)_{2m}$, $2m > 2$, $\sigma$ simple. An inflation $\pi$ of $\sigma$ is again in $I(321)$ if and only if $\pi$ is obtained through $\alpha_{i}=\alpha^{-1}_{\sigma^{-1}(i)}=\alpha^{-1}_{\sigma(i)}\,= 123\ldots m.$}\\

\begin{proof} By inflating a transposition of $\sigma$ through a substitution presenting an inversion,
we would obtain an inflation presenting a descending sequence of length 3, so not belonging to $Av(321)$.
\end{proof}

As a consequence of what is shown in 2.1, 2.2, 2.3, 2.6, we can affirm that the set $I(321)$, which is not strictly a {\it class} of permutations,  is a{n} almost \textit{ query-complete property}, in the sense of \cite{B}, pg.425, {Section} 1.
A property $P$ is said to be \textit{query-complete} if, for each simple permutation $\sigma \in P$, $\sigma$ of length $m$, there is a procedure to determine whether $\sigma[\alpha_{1},\ldots,\alpha_{m}]$ satisfies $P$ which requires only to know if each $\alpha_{i}$ satisfies $P$. The almost query-completeness allows to write the relations among generating functions described in the following {Section}.\\

\section{ Generating functions.}

We consider the following generating functions of subsets of $ I(321)$:\\
$\alpha$ the generating function of the involutions in $I(321)$ of the type 12;\\
$\beta$ the generating function of the involutions in $I(321)$ of the type 21;\\
$\gamma$ the generating function of simple involutions in $I(321)$, different from 1, 12 e 21;\\
$\delta$ the generating function of the involutions in $I(321)$ which are inflation of simple involution, of length  $n > 2$. \\Let finally $f$ be the generating function of { the whole set} $I(321)$.\\

On the basis of the structure theorems  and the properties of involutions recalled in {Section} 1, and of the properties of involutions avoiding 321 demonstrated in {Section} 2, adapting to our case the ideas exposed in \cite{B} e \cite{AA}, we write the following relations (1) for the generating functions:\\
 $$ \begin{cases} f = x \,+\, \alpha+\beta+\gamma+\delta\\
\beta = \frac{1}{1-x^{2}}-1\\
\alpha = (x\,+\,\beta+\gamma+\delta)(x\,+\,\alpha+\beta+\gamma+\delta)\end{cases}.\eqno{(1)}$$\\
The first equality reflects the structure theorem: an involution $\pi\neq 1$ is of the type 12, so enumerated by $\alpha$, or of the type 21, so enumerate by $\beta$, or it is simple or inflation of a simple, then enumerated by $\gamma$ or by $\delta$. {There exists only the involution $1$ of length $1$, hence in $f$ there is the summand $x$}. \\
{ The second equality follows from Proposition 2.2.}\\
The third relation is also based on a structure theorem: if $\pi$ is a type 12 involution, it can be uniquely expressed as $\pi=12[\sigma_{1},\,\sigma_{2}], \,\sigma_{1},\,\sigma_{2}\in I(321)$, with $\sigma_{1}$ not of type 12.\\

{Consider} the polynomial ring over the field $\cal K$, algebraic closure of the field of the rational functions. Let us try to determine the generating functions  $f,\alpha,\beta,\gamma,\delta$ starting from the polynomial  relations (1) (and we note that $x$ is in $\cal K$). We can say that if in the ideal generated by the relations there is a polynomial in $f$ only, { i.e. a relation between $f$ and $x$,} with coefficients in $\cal K$, then $f$ is an algebraic function.\\
The relations (1) do not suffice to determine a polynomial in $f$ only, so demonstrating the algebraicity of $f$, to which we arrive through the study of another property of the involutions avoiding 321.\\
\\Instead, if we want to calculate the generating function $\varphi$ of a subset $ I' \subset I(321)$, containing the only simple involutions 1, 12, 21, the relations are sufficient to obtain a generating function.\\\\
\textbf{Example 3.1.}  The relations (1) become (2):
$$ \begin{cases} \varphi= x+\alpha+\beta\\
\alpha= (x+\beta)(x+\alpha+\beta)\\
\beta = \frac{1}{1-x^{2}}-1\\
 \end{cases}.\eqno{(2)}$$
An easy calculation allows to express the rational generating function of $ I'$:
$$\varphi\,= \,\frac{x + x^{2} - x^{3}}{1 - x - 2 x^{2} + x^{3}},$$
whose expansion gives the coefficients 1, 2, 3, 6, 10, 19, 33, 61, 108, 197, 352, 638, 1145, 2069, 3721,6714, 12087
     $\ldots$ \\
It is easily found the recoursive formula for the coefficients of $\varphi$:
$$c_{n+3}= c_{n+2}+2c_{n+1}-c_{n},$$ with $c_{0}=1,\, c_{1}=2,\, c_{2}=3.$ ({S}ee Sloane, \cite {S}, A028495).\\

Back to the study of $I(321)$ and of new equations in order to express the generating function $f$, we demonstrate for the involutions in $I(321)$ the following properties, which connect the involutions of even length to the ones of odd length.\\
\\
{\bf Proposition 3.2}  \textit{Let $\pi'_{2m}\in I(321)_{2m}$, $m \geq 1$, con $\pi'_{2m}(1) = 1$. Then  $$\pi'_{2m}\, = \,12[1,\pi_{2m-1}]\quad\quad and \quad \quad  \pi_{2m-1}\in I(321)_{2m-1}.$$}
The property follows immediately from the structure theorems and Proposition 2.1.\\

We have already noticed that in $I(321))$ the involutions of type 21, the simple involutions of length greater than 2 and their inflations have no fixed points. So, such an involution $\pi$, $\pi\in I(321)_{2m}$, can be written as a product of  $m$ transpositions:
$$\pi\,=\,(m_1,M_1)(m_2,M_2)\ldots(m_m,M_m)\,$$
 where, for each $i=1,\ldots,m$,  $m_i< M_i$, and moreover $1=m_1<m_2<\cdots <m_m$ e $M_1<M_2<\cdots <M_m$.\\

>From Proposition 2.3. it {follows immediately} \\
\\
{\bf Proposition 3.3} \textit{For each involution $\pi\in I(321)_{2m}$, with $$\pi\,=\,(m_1,M_1)(m_2,M_2)\ldots(m_m,M_m),$$ another involution $\pi'$ can be defined, as below:
$$\pi'\,=\,(1)(m_2,M_1)(m_3,M_2)\ldots(m_m,M_{m-1})(M_{m}),$$
with $\pi'\in I(321)_{2m}$.}\\

(For example, for $\pi\,=\,\textbf{35}1\, \textbf{6}24\,= \,(1,3)(2,5)(4,6)\,\in I(321)_6$ {pro\-duct of three transpositions,} one has $\pi'\,=\,(1) (2,3)(4,5)(6) = 1\,\textbf{3}2\,\textbf{5}4\,\textbf{6}.$)\\
\\
{\textbf{Lemma 3.4} \textit{There is a bijection between the involutions $\pi \in I(321)_{2m}$, such that $\pi(1)\neq 1$, and the involutions such that $\pi(1)\,=\,1$.}\\

\begin{proof} The hypothesis $\pi(1)\neq 1$ says $M_1\neq 1$. So, such a $\pi \in I(321)_{2m}$ is either of type 21, or simple, or inflation of a simple  one of length {greater than 2}, or of type 12. In the first three cases, $\pi$ has no fixed points and Proposition 3.3 provides an involution $\pi' \in I(321)_{2m}$ such that $\pi' (1)=1$. If $\pi$ is of type 12, it has the canonical decomposition of Proposition 1.3, $\pi = 12[\alpha_{1},\alpha_{2}]$, with $\alpha_{_1}$ not of type 12. From the previous propositions $\alpha _1$ has even length and no fixed points, then $\pi^\prime =12[\alpha ^\prime_1, \alpha _2]$, where $\alpha ^\prime_1$ is obtained from $\alpha_1 $ by means of Proposition 3.3, is again an involution of $I(321)_{2m}$ with $\pi^\prime(1)=1$. \\
Viceversa, let  $\pi'\,\in I(321)_{2m} $ be such that  $\pi'(1)=1$: then $\pi'$ has at least a second fixed point, (recall that $\pi'_{2m} $ must have an even number of fixed points). Let $M_{i}\neq 1$ be the first of the fixed points following 1: then $\pi'\,=$ $(1,1)(m_2,M_2)(m_3,M_3)\cdots (m_{i-1},M_{i-1})(M_{i},M_i)(m_{i+1},M_{i+1})\ldots (m_h,M_{h}).$\\
The involution  $\pi'$ is precisely obtained from \\
\centerline{$\pi\,=\,(1,M_2)(m_2,M_3)\ldots (m_{i-1},M_{i})\ldots (m_h,M_{h}).$}
Indeed if $i<h$, then $\pi=\,12[\alpha _1,\alpha_2]{,}$ with\\ \centerline{$\alpha_1=\,(1,1)(m_2,M_2)(m_3,M_3)\ldots$ $(m_{i-1},M_{i-1})(M_{i},M_{1}){,}$}
{so} $\,12[\alpha^\prime_1,\alpha_2]\,$ is just $\pi^\prime$;\\ if $i\,=\,h,$ $\,\pi\,$ is not of type 12 and $\,\pi^\prime\,$ is obtained from $\,\pi\,$ by the construction of Proposition 3.3.
\end {proof}
\\
\textbf{Example } Let $\pi = \textbf{4\,6\,7}\,1\, \textbf{8}\,2\,3\,5 \,(\textbf{10})\,9\;=\;12[\,46718235,\,21\,]\,=\,(1,4)(2,6)(3,7)$ $(5,8)(9,(10))\,$, where $\textbf{4\,6\,7}\,1\, \textbf{8}\,2\,3\,5\;$ is expansion of the simple involution  $\textbf{3\,5}\,1\, \textbf{6}\,2\,4\,$.\\Then $\pi^\prime= \;12[\,14627358,\,21\,]\,$ $1\,\textbf{4\,6}\,2\,\textbf{7}\,3\,5\,\textbf{8\,(10)}\,9\,=\,(1)(2,4)(3,6)(5,7)(8)(9(10))$.\\

>From Proposition 3.2 and Lemma 3.4  we finally derive\\
\\
\textbf{Theorem 3.5} \textit{For the orders of the sets $I(321)_{2m}$ and $I(321)_{2m-1}$, the equality holds $$|I(321)_{2m}|=2|I(321)_{2m-1}|.$$}

We now introduce the generating functions $\varepsilon$ and $\omega$, respectively of the sets of the involutions of even length  and of the involutions of odd length.\\
\\
\textbf{Theorem 3.6}. \textit{Let $f$, $\varepsilon$, $\omega$ be the generating functions respectively of $I(321)$, of $\bigcup I(321)_{2m}$ and of $\bigcup I(321)_{2m+1},\,m\in N$. We have the relations
$${\varepsilon }+ {\omega} = f  \, ,\qquad { \varepsilon} = 2x {\omega}.$$ }

 \section{Algebraicity of the generating functions of $I(321)$ and of some of its subsets.}

>From Propositions 2.1, 2.2, 2.3 we derive that in $I(321)$ the involutions of type 21, the simple ones and their inflations have all even length, while  an involution  of odd length can only be 1 or of type 12. Then we can write the following relations (3) between the generating functions $\varepsilon$ e $\omega$:
$$\,{\omega}\,=\,x\,+\,x{\varepsilon }\,+(\beta+\gamma +\delta){\omega };\quad \quad\varepsilon=\beta+\gamma+\delta\,+\, (\beta+\gamma+\delta){\varepsilon }\,+\,x{\omega}.\eqno{(3)}$$
The first equality describes the property that the involutions of odd length avoiding 321 are of the following kinds: \\$\pi_{1} = 1$,\\ $\pi_{2m+1} = 12[1,\pi_{2m}]$, con $\pi_{2m}\in I(321)$, so enumerated by $x{\varepsilon }$,\\ $\pi'_{2m+1} = 12[\alpha_{1},\alpha_{2}]$, where $\alpha_{1}$ has even length and is not of type 12, and $\alpha_{2} \in I(321)$ has odd length.\\
 The second equality describes the  property already recalled that the involutions in $I(321)$ of type 21, the simple involutions of length greater than $2$ and their inflations have an even length; moreover the following involutions have an even length:\\
 $\pi_{2m+2} = 12[\alpha_{1},\alpha_{2}]$, where $\alpha_{1}$ has even length and is not of type 12, and $\alpha_{2}\in  I(321)$ has  even length;\\
$\pi'_{2m+2} = 12[1,\pi_{2m+1}]$,  $\pi_{2m+1}\in I(321)$, so enumerated by $x{\omega}$.\\
Equalities (1) and (3), and the property ${\varepsilon} = 2x {\omega}$ give the system (4):
$$ \begin{cases} f = x \,+\, \alpha+\beta+\gamma+\delta\\
\beta = \frac{1}{1-x^{2}}-1\\
\alpha = (x\,+\beta+\gamma+\delta)(x\,+\,\alpha+\beta+\gamma+\delta)\\
f={\omega}+{\varepsilon}\\
{\omega}\,=\,x\,+\,x{\varepsilon }\,+(\beta+\gamma +\delta){\omega }\\
{\varepsilon }=\beta+\gamma+\delta +\,(\beta+\gamma+\delta){\varepsilon }\,+\,x{\omega}\\
{\varepsilon }=2x{\omega }\end{cases}.\eqno{(4)}$$
The relations (4) allow to find the polynomial in $f$ and $x$,

$$-f + x - f^{2} x + 2 x^{2} + 4 f x^{2} +
        2f^{2} x^{2},$$
which leads to the generating function
$$f=\frac{1 - 4 x^{2} -\sqrt{1 - 4 x^{2}}}{2(-x + 2x^{2})}
,$$
whose expansion gives the coefficients\\
\centerline {$1, 2, 3, 6, 10, 20, 35, 70, 126, 252, 462, 924, 1716, 3432, 6435,
      12870, 24310,\ldots$}\\
({wich are the } central binomial coefficients, see \cite {S}, A001405).\\
{ In this way we {obtain} a new proof of\\
\\
\textbf{Theorem 4.1} ({S}ee \cite{ss}, \cite{M}.)  \textit{The size of $I(321)_n$ is $ n \choose \lfloor n/2 \rfloor $.}\\\\
{This equality was firstly proved in \cite{ss} by means of the Standard Young Tableaux; a combinatorial proof based on a class of Dyck paths was presented in \cite {M}.}\\

>From (4) one also derives for $\varepsilon$ $${\varepsilon }\,=\, \frac{-1-4x^{2}-\sqrt{1-4x^{2}}}{-1+4x^{2}}.$$

{
{We note that in} \cite{B} and \cite{AA} it is shown that a permutation class  with only \textit{finitely many} simple permutations has a  readily computable algebraic generating function (besides more general results, always in the ipothesis of the existence of finitely many simple permutations). But the set $I(321)$  has \textit{infinitely many }simple involutions, as noted from the beginning ( {in} Proposition 2.5.) This {situation} gives another example of a {case} where the generating function is algebraic, whereas in presence of infinitely many simple permutations.\\

The system (4) allows to express also the polynomial in $\alpha$ only,
$$ \alpha + \alpha x + \alpha ^2 x - x^2 - 4 \alpha x^2 - 2 \alpha ^2 x^2 - 4x^3 -
        4 \alpha x^3 - 4 x^4 , $$
which leads to
 $$\alpha = \frac{1 + x - 4 x^2 - 4 x^3 -\sqrt{1 + 2 x - 7 x^2 - 12x^3 + 16 x^4 + 16 x^5 - 16 x^6}}{2(-x + 2 x^2)}$$
 whose expansion gives the coefficients\\
1, 3, 5, 10, 18, 35, 65, 126, 238, 462, 882, 1716, 3300, 6435,
      12441, 24310, 47190, 92378, 179894 $\ldots$ ({see \cite {S}, A107232}). \\

The use of the Gr\"obner Basis' theory not only is useful for the calculations, but also allows to affirm that {in the ideal generated by}  (4) there exists no polynomial in $\,\gamma\,$ only, nor in $\,\delta\,$ only.  While there exists a polynomial in  $\zeta\,=\,\gamma+\delta$ only, i.e. a relation between $\gamma +\delta$ and $x$,
$$(-1+4x^{2}-3x^{4})(\zeta)+(1-2x^{2}+x^{4})(\zeta )^{2}+x^{6}\eqno{(5)}$$
which leads to $$\zeta\,=\,\gamma + \delta \,=\,\frac{1 - 4 x^{2} + 3 x^{4}-\sqrt{1 - 8 x^{2} + 22 x^{4 }- 28x^{6} + 17 x^{8} - 4 x^{10}}}{2(1 - 2 x^{2} + x^{4})}$$
 whose expansion gives the coefficients\\
 1, 0, 4, 0, 13, 0, 41, 0, 131, 0, 428, 0, 1429, 0, 4861. \\(The non zero terms are Catalan numbers, see \cite {S} , A001453).\\
We note that at this point one cannot { know if} $\gamma$ and $\delta$ are algebraic. \\

 \section{The generating functions of {the} simple involutions and their inflations.}

In order to calculate the generating function $\gamma$ of the simple involutions in $I(321)$ we adapt to the involutions' case the condiderations of \cite{a},  3.\\
{ Recall that for each of the simple involutions of length greater than $2$, enumerated by $\gamma $, its expansions are obtained, among the involutions enumerated by $\zeta$,  by inflating the pairs of elements of some transpositions by means of increasing sequences $\,12\,$, $\,123\,$, ... ,$\,1\cdots n\,$, ...  ({as seen in }Proposition 2.6).\\ Then $\zeta (x^2)=\gamma \left ( {{x^2}\over {1-x^2}}\right )$ and conversely $\gamma (x^2)= \zeta \left ( {{x^2}\over{1+x^2}}\right )$.}\\

 Thus the substitution in the polynomial (5) leads to the polynomial
 $$ (1+x^2)z^2+(-1+x^2+2x^4)z+x^6 \, ,$$
 hence to the algebraic generating function
$$ \gamma = {{1 - x^2 - 2x^4 -\sqrt{1 - 2 x^2 - 3 x^4}}\over{2(1 + x^2)}} \, , $$
{whose expansion gives the coefficients}\\
 \centerline {  $1, 0, 1, 0, 3, 0, 6, 0, 15, 0, 36, 0, 91, 0, 232, 0, 603, 0, 1585, \dots $ }\\{The non zero coefficients are Riordan's numbers, see \cite {S} , A005043.}\\

For the function $\delta$, generating the inflations of simple involutions, the  initial coefficients of the expansion are determined through the difference $\delta$ = $\zeta$ - $\gamma$, namely:\\
3, 0, 10, 0, 35, 0, 116, 0, 392, 0, 1338, 0, 4629, 0, 16192, 0,
      57200, 0, 203798, 0, 731601, 0, 2643902, 0, 9611747, 0, 35130194, 0,
      129018797, 0, 475907912, 0, 1762457594....\\
(The sequence is not listed in \cite {S}{)}.\\

Also the algebraic function $f\, -\, \gamma$, generating function of the set $$I(321)\,\cap\,Av(2413,\,3142),$$  (recall that  Av(2413, 3142) is the class of the  \textit{separable permutations}), leads to a sequence not listed in \cite {S}:\\
1, 2, 3, 6, 10, 19, 35, 69, 126, 249, 462, 918, 1716, 3417, 6435,
      12834, 24310, 48529, 92378, 184524 $\dots$\\

For the coefficients $\delta_{2n}$ of the generating function $\delta$ we give the following property.\\
\\
{\bf Theorem 5.1} {\it The coefficients $\delta_{2n}$ of  $\delta$ satisfy
$$\delta_{2n}\,=\,\sum_{i=1}^{n-3}\, \gamma_{2(n-i)}\sum_{j=0}^{i-1} {i-1\choose j}{n-i\choose j+1} .$$}
(The usual convention holds, that ${h\choose k}=0$ when $h<k$. Because, in the last addendum of the sum for $i$, one has $i= n-3$ and then $n-i=3$, the sum for $j$ has only  three addenda not zero; in the second from last only four, and so on.)\\
 \begin{proof} The involutions enumerated by $\delta_{2n}$ are the inflations of the simple involutions $\sigma_{m} \in I(321)$, of length $m,\,\,{\rm con}\,\, 6 \leq m < 2n$; from  Proposition 2.6 such an inflation is obtained through ascending sequences, the same for the terms of each transposition.\\
 Each one of the  $\gamma _{2(n-1)}$ simple involutions of length $2(n-1)$ is inflated to one of length $2n$ by substituting the sequence 12 in each of the terms of a transposition: having $\gamma _{2(n-1)}$ the number of $n-1= {n-1\choose 1}$ transpositions, the contribution of the simple permutations of length $\;2(n-1)\;$ is $\;\gamma _{2(n-1)} {n-1\choose 1}\;,$ which we write in the form $\gamma _{2(n-1)}{0\choose 0} {n-1\choose 1}$.\\

 For each of the $\gamma _{2(n-2)}$ simple involutions of length $2(n-2)$ we can obtain an involution of length $2n$ by sustituting the two terms of the $n-2\choose 1$ transpositions with the sequence 123, or by sustituting the four terms of a couple of transpositions (wich are $n-2\choose 2$) with the sequence 12. So the contribution of the simple involutions of length $\;2(n-2)\;$ is $\;\gamma _{2(n-2)}\left[ {n-2\choose 1}+{n-2\choose 2}\right ],$
 which we write in the form $\gamma _{2(n-2)}\left[{1\choose 0}{n-2\choose 1}+{1\choose 1}{n-2\choose 2}\right ]$.\\

The general formula follows, by observing that each one of the $\gamma _{2(n-i)}$ simple involutions of length ${2(n-i)}$, so containing $n-i$ traspositions, contributes by means of:\\
one transposition, whose terms are substituted by the ascending sequence $\,1\cdots i+1$; for each one of the $n-i\choose 1$ transpositions the number of inflations is $1={i-1\choose 0}$, with $j=0$;\\
a couple of transpositions, where the sequence $1\cdots h+1$ is substituted for the terms of the first couple, and the sequence $\;1\cdots i-h+1\;$ for the terms of the second; so, for each one of the  $n-i\choose 2$ couples we have $i-1={i-1\choose 1}$ different inflations, with $j=1$;\\
a tern of transpositions, where the sequence $\;1\cdots h+1\;$ is substituted for the terms of the first, $\;1\cdots k+1\;$ for the terms of the second, $\;1\cdots i-(h+k)+1\;$  for the terms of the third; so, for each one of the $n-i\choose 3$ terns the different inflations are $i-1\choose 2$ (with $j=3-1=2$), namely as many as the couples of numbers $h,\,k$ chosen between $1$ e $i-1$.\\
And so { on}. \end{proof}

\section{The graphic of the simple involutions.}

Consider the graphic of an involution $\pi \in I(321)$ with no fixed points: the sequence of maxima is represented over the line $y=x$, the sequence of minima under the line. Let's connect the points of the graph in the order their ordinates possess in the permutation $\pi$. We call \textit{plot of the involution} the drawing so obtained. \\Define two maxima (or two minima) to be \textit{up-connected} (respectively \textit{down-connected}) when connected through a step of the drawing neither crossing the line $y=x$ nor containing another maximum (or minimum), therefore consecu{-}tive in the permutation. In this case, we also say that the plot has an \textit{upper connection} (or a \textit{lower connection}).\\\

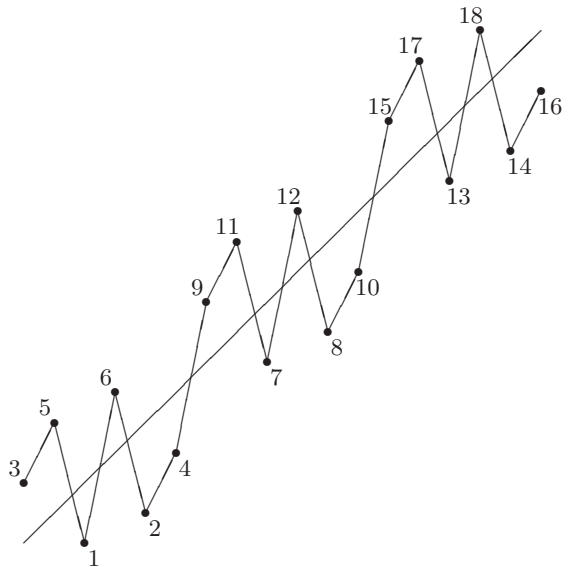
\begin{figure}[h]

\begin{center}
 \setlength{\unitlength}{4mm}
 \begin{picture}(30,18)

\put(6,0){
\begin {picture}(0,0)
\put(0,2){\circle*{0.3}} \put(1,4){\circle*{0.3}} \put(2,0){\circle*{0.3}}\put(3,5){\circle*{0.3}} \put(4,1){\circle*{0.3}} \put(5,3){\circle*{0.3}}
\put(6,8){\circle*{0.3}} \put(7,10){\circle*{0.3}} \put(8,6){\circle*{0.3}}\put(9,11){\circle*{0.3}} \put(10,7){\circle*{0.3}} \put(11,9){\circle*{0.3}}
\put(12,14){\circle*{0.3}} \put(13,16){\circle*{0.3}} \put(14,12){\circle*{0.3}}\put(15,17){\circle*{0.3}} \put(16,13){\circle*{0.3}} \put(17,15){\circle*{0.3}}

\put(0,2){\line(1,2){1}} \put(1,4){\line(1,-4){1}} \put(2,0){\line(1,5){1}} \put(3,5){\line(1,-4){1}} \put(4,1){\line(1,2){1}}
\put(6,8){\line(1,2){1}} \put(7,10){\line(1,-4){1}} \put(8,6){\line(1,5){1}} \put(9,11){\line(1,-4){1}} \put(10,7){\line(1,2){1}}
\put(12,14){\line(1,2){1}} \put(13,16){\line(1,-4){1}} \put(14,12){\line(1,5){1}} \put(15,17){\line(1,-4){1}} \put(16,13){\line(1,2){1}}
\put(5,3){\line(1,5){1}} \put(11,9){\line(1,5){1}}

\put(0,0){\line(1,1){17}}

\put(-0.3,2.5){\makebox(0,0){\small $3$}}
\put(0.7,4.5){\makebox(0,0){\small $5$}}
\put(2.3,-0.5){\makebox(0,0){\small $1$}}
\put(2.7,5.5){\makebox(0,0){\small $6$}}
\put(4.3,0.5){\makebox(0,0){\small $2$}}
\put(5.3,2.5){\makebox(0,0){\small $4$}}

\put(5.7,8.5){\makebox(0,0){\small $9$}}
\put(6.7,10.5){\makebox(0,0){\small $11$}}
\put(8.7,11.5){\makebox(0,0){\small $12$}}

\put(11.7,14.5){\makebox(0,0){\small $15$}}
\put(12.7,16.5){\makebox(0,0){\small $17$}}
\put(14.7,17.5){\makebox(0,0){\small $18$}}

\put(8.3,5.5){\makebox(0,0){\small $7$}}
\put(10.3,6.5){\makebox(0,0){\small $8$}}
\put(11.3,8.5){\makebox(0,0){\small $10$}}

\put(14.3,11.5){\makebox(0,0){\small $13$}}
\put(16.3,12.5){\makebox(0,0){\small $14$}}
\put(17.3,14.5){\makebox(0,0){\small $16$}}

\end{picture}}

\end{picture}
\end{center}

\caption{Plot of $\pi=3516249(11)7(12)8(10)(15)(17)(13)(18)(14)(16)$}\label{Figure S4}
\end{figure}

\noindent The involution of Figure 1 is not simple, indeed $\pi =12[351624\,,\,12[351624\,,\,351624]]$; the pairs of points $3,5$; $9,11$; $15,17$ are up connected; the pairs $2,4$; $8,10$; $14,16$ are down connected.\\
\\
\textbf{Proposition 6.1.} \textit{Let $\pi \in I(321),\,\pi \, = \, (m_1,M_1)(m_2,M_2)\cdots (m_m,M_m)$.  If two minima, $m_i,\, m_{i+1}$, are down-connected, then the corresponding maxima, $M_i,\, M_{i+1}$, are consecutive integers. \\
Conversely, if two maxima, $M_i,\, M_{i+1}$, are consecutive integers, then the corresponding minima, $m_i,\, m_{i+1}$, are down-connected.\\
(The same holds with maximum  replaced by minimum, and up by down.)}\\

\begin{proof} If $M_i +1\neq M_{i+1}$, another integer $M$ would exist, with $M_i<M< M_{i+1}$, wich can be neither a maximum nor a minimum.\\
Conversely, if $M_i +1 = M_{i+1}$, then $m_i$ and $m_{i+1}$ are down-connected, because otherwise another maximum  $M$ should exist, with $M_i +1\,<\,M\,<\,M_{i+1}$.\end{proof}

It is immediately obtained by Proposition 6.1 the following property of the simple involutions of $I(321)$:
\\\\
\textbf{Proposition 6.2.} \textit{Let $\sigma \in I(321)$, $\sigma$ simple. Then the plot of $\sigma$ has no couples of upper and lower connections symmetric with respect to $y=x$, therefore if two maxima are up-connected, the corresponding minima are not down-connected.} \\

 The converse is not true. However we have:
\\\\
\textbf{Theorem 6.3.} \textit{If the plot of an involution $\pi\in I(321)$ has no couples of symmetric connections, then {either} $\pi$ is  a simple involution, or $\pi$ is an involution of type 12, with $\pi\,=\,12[\alpha_{1},\,\alpha_{2}],$ where $\alpha_{1}$ is simple.}\\

\begin{proof} { Recalling the first equation of (1), if} $\pi$ were of type 21, by Proposition 2.1 its plot would have at least a couple of symmetric upper and lower connections. If $\pi$ were an expansion of a simple involution $\sigma \neq 12$, by Proposition 2.6 its plot again would have at least a couple of symmetric upper and lower connections (deriving by the inflation of  a transposition). \\{Then either $\pi$ is simple or it is of type 12. In the last case $\pi\,=\,12[\alpha_{1},\,\alpha_{2}]$ where $\alpha_{1}$ must be simple, while $\alpha _2$, satisfying the hypothesis of the theorem, must be again either simple or of type 12}.\end{proof}

\section{Simple involutions in $I(321)$ and {short} Motzkin paths with no horizontal steps at level $0$.}

Consider the connection between the coefficients of the expansion of the function $\gamma$ and the  Riordan's numbers {also called Motzkin sums}, through the interpretation presented in \cite {S}, given by Emeric Deutsch (2003): \textit{the coefficients \centerline{ $\,\,1, 0, 1, 1, 3,  6,  15,  36,  91,  232,  603, 1585, \ldots$} \\enumerate  the Motzkin paths of length $n$ with no horizontal steps at level $0$.}\\
We call them conventionally {\it short}  Motzkin paths because, in the following {Section}, we shall instead refer to Motzkin paths of length $2n+2$.\\
Such an interpretation, together with Proposition 6.2, leads to the complete description of the plot of simple involutions of $I(321)$, through the following procedure, that defines a bijection between the simple involutions of $I(321)_{2n+2}$ and the Motzkin paths with no horizontal steps at level $0$ of length $n$.\\

Written $\sigma_{2n+2} $ in the form $\sigma\,=\,(M_1,m_1)(M_2,m_2)\ldots(M_{n+1},m_{n+1})$, one can consider the  sequence of $n+1$ odd integers  $\left\{1,3,\dots,1 \right\}$, that  describes the number of times the plot of $\sigma_{2n+2} $ crosses the line $y=x$ between each maximum and the respective minimum. For instance, one has:\\
for $\sigma_{6}\,=\,\textbf{35}1\textbf{6}24 = (31)(52)(64)$ the corresponding sequence is $\left\{1,3,1 \right\}$; \\
for $\sigma_{8}\,=\,\textbf{35}1\textbf{7}2\textbf{8}46 = (31)(52)(74)(86)$ the corresponding sequence is $\left\{1,3,3,1 \right\}$; in the appendix we list simple involutions and their sequences for $n=10,12,14$.
\\\\
{{\bf Proposition 7.1.} \textit{Let $\sigma_{2n+2}\in I(321)_{2n+2} $ with $\{s_1,\ldots ,s_{n+1}\}$ the associated sequence. Then\\i) $s_i$ is odd (for $i=1,\ldots,n+1$);\\ ii) $s_1=s_{n+1}=1$;\\iii) if $s_i=1$ then both $s_{i-1}$ and $s_{i+1}$ are different from $1$; \\
iv) $|s_{i+1}-s_i|\leq 2$ (for $i=1,\ldots,n$).}\\

\begin{proof} All of the properties derive from the structure of the simple involutions and from Proposition 6.2, almost immediately for \textit{i),ii),iii)}, remembering that in a simple involution no maximum  can be ad{j}acent to its minimum. \\As for \textit{iv)},  it is to  be noted that going from $M_{i }$ to $m_{i }$ with ${s}_{i}$ crossings always means that $M_{i+1 }$ is followed by $m_{i }$. Now, one has:\\either $m_{i }$ and $m_{i+1 }$ are connected, so ${s}_{i+1}={s}_{i} -2$, because in this case $M_{i }$ and $M_{i+1 }$ cannot be connected;\\
or $m_{i }$ and $m_{i+1 }$ are not connected, and then: either $M_{i }$ and $M_{i+1 }$ are not connected, so ${s}_{i }$ = ${s}_{i+1 }$, or $M_{i }$ and $M_{i+1 }$ are  connected, so ${s}_{i +1} = {s}_{i }+2$.
\end{proof}
}
\\\\
We call \textit{admissible sequence} for a simple involution a sequence of consecutive or repeated odd numbers satisfying the claims of Proposition 7.1.\\ To each admissible sequence $\left\{s_{i}\right\}$ of length  $n+1$ we associate the Motzkin path of length $n$ presenting up, down or horizontal steps  depending on $s_{i}< s_{i+1}$, $s_{i}>s_{i+1}$, $s_{i}=s_{i+1}$. \quad\\
Conversely, to each Motzkin path of length $n$ with no horizontal steps at level $0$ we associate the  admissible sequence of odd numbers, of length  $n+1$, $\left\{s_{i}\right\}$ = $\left\{1,3,\dots,1 \right\}$, where  $s_1=1$, $s_{i+1} = s_{i}+2$ if  the Motzkin path has an up step,
$s_{i+1} = s_{i}-2$ if the Motzkin path has a down step, $s_{i+1} = s_{i}$ if the Motzkin path is horizontal.

\begin{figure}[h]

\begin{center}
 \setlength{\unitlength}{4mm}
 \begin{picture}(30,7)

\put(3,0){
\begin {picture}(0,0)
\put(0,2){\circle*{0.3}} \put(1,4){\circle*{0.3}} \put(2,0){\circle*{0.3}}\put(3,5){\circle*{0.3}} \put(4,1){\circle*{0.3}} \put(5,3){\circle*{0.3}}
\put(0,2){\line(1,2){1}} \put(1,4){\line(1,-4){1}} \put(2,0){\line(1,5){1}} \put(3,5){\line(1,-4){1}} \put(4,1){\line(1,2){1}}
\put(0,0){\line(1,1){5}}

\put(-0.3,2.5){\makebox(0,0){\small $3$}}
\put(0.7,4.5){\makebox(0,0){\small $5$}}
\put(2.3,-0.5){\makebox(0,0){\small $1$}}
\put(2.7,5.5){\makebox(0,0){\small $6$}}
\put(4.3,0.5){\makebox(0,0){\small $2$}}
\put(5.3,2.5){\makebox(0,0){\small $4$}}

\end{picture}}

\put(15,3){\makebox(0,0){\small $\{131\}$}}

 \put(22,2){
 \begin{picture}(0,0)
 \put(0,0){\circle*{0.3}} \put(2,2){\circle*{0.3}} \put(4,0){\circle*{0.3}}
\put(0,0){\line(1,1){2}} \put(2,2){\line(1,-1){2}}
 \end{picture}}
\end{picture}
\end{center}

\caption{The only simple involution $\,351624\,= (1,3)(2,5)(4,6)\in I(321)_6$ and the corresponding Motzkin path of length 2.}\label{Figure 2}
\end{figure}
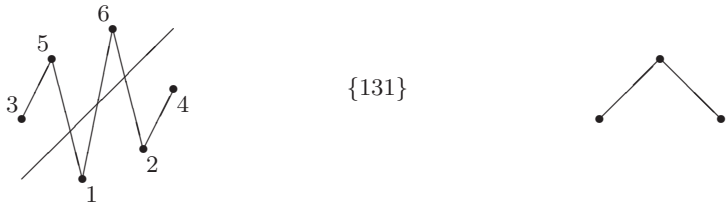

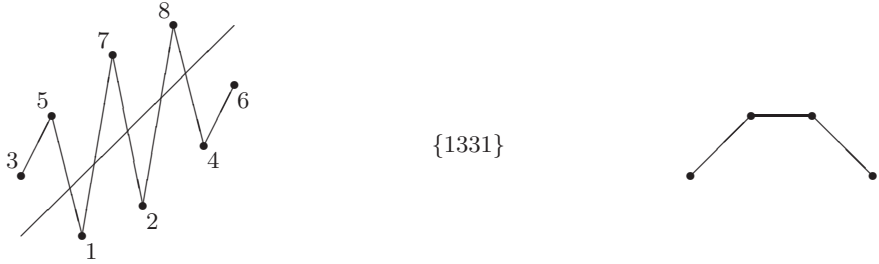
\begin{figure}[h]

\begin{center}
 \setlength{\unitlength}{4mm}
 \begin{picture}(30,8)

\put(0,0){
\begin {picture}(0,0)
\put(0,2){\circle*{0.3}} \put(1,4){\circle*{0.3}} \put(2,0){\circle*{0.3}}\put(3,6){\circle*{0.3}} \put(4,1){\circle*{0.3}} \put(5,7){\circle*{0.3}} \put(6,3){\circle*{0.3}} \put(7,5){\circle*{0.3}}
\put(0,2){\line(1,2){1}} \put(1,4){\line(1,-4){1}} \put(2,0){\line(1,6){1}} \put(3,6){\line(1,-5){1}} \put(4,1){\line(1,6){1}}
\put(5,7){\line(1,-4){1}} \put(6,3){\line(1,2){1}}
\put(0,0){\line(1,1){7}}

\put(-0.3,2.5){\makebox(0,0){\small $3$}}
\put(0.7,4.5){\makebox(0,0){\small $5$}}
\put(2.3,-0.5){\makebox(0,0){\small $1$}}
\put(2.7,6.5){\makebox(0,0){\small $7$}}
\put(4.3,0.5){\makebox(0,0){\small $2$}}
\put(4.7,7.5){\makebox(0,0){\small $8$}}
\put(6.3,2.5){\makebox(0,0){\small $4$}}
\put(7.3,4.5){\makebox(0,0){\small $6$}}

\end{picture}}

\put(15,3){\makebox(0,0){\small $\{1331\}$}}

 \put(22,2){
 \begin{picture}(0,0)
 \put(0,0){\circle*{0.3}} \put(2,2){\circle*{0.3}} \put(4,2){\circle*{0.3}} \put(6,0){\circle*{0.3}}
\put(0,0){\line(1,1){2}} \put(2,2){\line(1,0){2}} \put(4,2){\line(1,-1){2}}

 \end{picture}}
\end{picture}
\end{center}

\caption{The only simple involution $\;35172846\,\in I(321)_8$ and the corresponding Motzkin path of length 3.}\label{Figure 3}
\end{figure}

To each Motzkin path of length $n$ with no horizontal steps at level $0$ we associate the plot of a simple involution $\sigma_{n+2}$ and the involution $\sigma_{n+2}$ itself, univocally determined in the following way.\\
Define \textit{ local maximum point} of the Motzkin path a point $P_{i}$ such that $P_{i-1}\leq P_{i}$ and $P_{i+1}\leq P_{i}$.\\
Starting from the left, consider the number $N_{1}$ of the ascending steps to the left of the first  local maximum point $P_{1}$. Then $N_{1}+1$ determines the number of the maxima preceding $m_{1}=1$ in the involution. Similarly, if $N_{2}$ is the number of the descending steps at the right of the last  local maximum point, $N_{2}+1$ determines the number of the minima following $M_{n+1}\,=\,2n+2$ in the involution. Between the first and the last  local maximum point, each up step determines an up connection between two maxima of the involution, while each down step determines a down connection between two minima. If the Motzkin path has an only maximum, the plot has no up or down connections different from the ones regarding the first maxima and the last minima (for example, in the case $\sigma_{14}=\textbf{579(11)}1 \textbf{(12)}2 \textbf{(13)}3 \textbf{(14)}468(10)$).\\
Through the sequence $\left\{1,3,\dots,1 \right\}$ corresponding to the plot, one easily sees that a bijection yields,
giving a combinatorial interpretation for the cardinali{-}ty's equality  of the two considered sets.
\\ \\Connecting Theorem 6.3 and this combinatorial interpretation, we obtain \\\\
\textbf{Theorem 7.2} \textit{The plot of an involution $\pi_{n+2}\in I(321)$, presenting no couples of symmetric upper and lower connections, corresponds to a simple involution if and only if generates a  Motzkin path of length $n$ with no horizontal steps at level $0$.}\\

\begin{proof} By Theorem 6.3, or $\pi$ is a simple involution, so generating such a Motzkin path, or is of type 12, with $\pi\,=\,12[\sigma_{1},\,\sigma_{2}],$ where $\sigma_{1}$ is simple. In the second case, by the natural generalization of the procedure introduced for the simple involutions, $\pi$  generates a  Motzkin path, with horizontal steps at level $0$  separating the simple components of the permutation. (See example  in fig. 4).\\
In such a way we {can immediately recognize in the path the simple components of the involution.}
\end{proof}

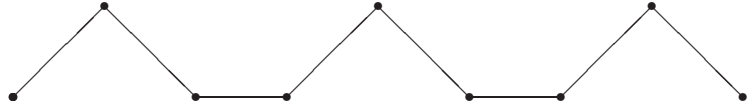
\begin{figure}[h]

\begin{center}
 \setlength{\unitlength}{4mm}
 \begin{picture}(30,6)

 \put(15,5){\makebox (0,0) {$\{131131131\}$}}
 \put(4,0){\begin{picture} (0,0)
 \put(0,0){\circle*{0.3}} \put(3,3){\circle*{0.3}} \put(6,0){\circle*{0.3}} \put(9,0){\circle*{0.3}} \put(12,3){\circle*{0.3}} \put(15,0){\circle*{0.3}} \put(18,0){\circle*{0.3}} \put(21,3){\circle*{0.3}} \put(24,0){\circle*{0.3}} \put(0,0){\circle*{0.3}} \put(0,0){\circle*{0.3}} \put(0,0){\circle*{0.3}}

\put(0,0){\line(1,1){3}} \put(3,3){\line(1,-1){3}} \put(6,0){\line(1,0){3}} \put(9,0){\line(1,1){3}} \put(12,3){\line(1,-1){3}} \put(15,0){\line(1,0){3}} \put(18,0){\line(1,1){3}} \put(21,3){\line(1,-1){3}}
 \end{picture}}
 \end{picture}
 \end{center}
 \caption{decomposition of $\,351624\,9(11)7(12)8(10)\,(15)(17)(13)(18)(14)(16)\;$}\label{Figure 4}
 \end{figure}

{ \textbf{Remark 1.} Given $\sigma_{2n+2} \in I(321)_{2n+2},$ it is of interest the problem of deciding which simple involutions are contained in  $\sigma_{2n+2}$ as patterns.  Considering the sequence $s$ = $\left\{1,3,\dots,1 \right\}$, of length $n+1$, associated to  a simple involution $\sigma_{2n+2} \in I(321)_{2n+2}$, one can actually determine the simple involutions  contained in $\sigma_{2n+2}.$ Those are precisely the simple involutions corresponding to admissible sequences,  obtained as
admissible subsequences of $s$,  of length $< n+1$.\\

For instance, $\sigma_{12}=\textbf{3 5}1 \textbf{7}2 \textbf{9}4 (\textbf{11})6 (\textbf{12})8(10)$, whose corresponding sequence is $\left\{133331\right\}$, contains the simple involutions of inferior length:\\
 $\textbf{35}1\textbf{7}2 \textbf{9}4 (\textbf{10}) 68,$ corresponding to $\left\{13331\right\}$, containing $\sigma_{8}$ e $\sigma_{6}$, respectively corresponding to $\left\{1331\right\}$ and $\left\{131\right\}$.\\

On the contrary,  $\sigma_{12}=\textbf{3 5}1 \textbf{7}2 \textbf{9}4 (\textbf{11})6 (\textbf{12})8(10)$ does not contain the involutions \textbf{36}1\textbf{79}24 (\textbf{10})58, corresponding to
$\left\{13131\right\}$, and $\textbf{468}1\textbf{9}2 (\textbf{10})357,$ corresponding to the sequence $\left\{13531\right\}$, not having the requested property.}\\

\textbf{Remark 2.} Analogously, starting from $\sigma_{2n+2} \in I(321)_{2n+2}$ and its associated sequence $ s$, one can build all the simple involutions $\sigma_{2n+4} \in I(321)_{2n+4}$, containing
as pattern the involution $\sigma_{2n+2}$ .\\
\section{Simple involutions in $I(321)$ and Dyck paths.}

In this {Section} we show how through the use of  Dyck paths another procedure can be obtained to determine whether an involution is simple, and to calculate simple involutions.\\
In \cite{AAA}, Proposition 4, involutions in $I(321)$ are characterized in terms of \textit{labelled Motzkin paths}. How to associate an involution $\pi_{n}$ of length $n$ to a  labelled Motzkin path  $(M,\lambda)$, where $M$ denotes a Motzkin path of length  $n$ and $\lambda$ denotes a labelling of { its} down steps, it is discussed thoroughly in \cite{AAA}, {Section} 3, pg.3. We only recall here that, given an involution $\pi$ in the form $\pi\,=\,(m_1,M_1)(m_2,M_2)\ldots(m_m,M_m)\,$
 whith the $m_i$ written in increasing order,  a labelled Motzkin path is defined as follows. For every $i=1,\ldots,n,$ \\
- if $i$ is a fixed point for $\pi$, take a horizontal step in the  path;\\
- if $i$ is the first element af a transposition, take an up step in the   path; \\
- if $i$ is the second element af a transposition, take a down step in the   path, \textit{labelled} with $h$, if $i$ is in the $h$-th position among integers greater than or equal to $i$ in the cycle decomposition of $\pi.$ \\ The following  proposition holds, from which we derive the characterization in Theorem 8.2.\\

\textbf{Proposition 8.1.} (See \cite{AAA}, Proposition 4.) \textit{Let $\pi_{n}$ be an involution with $(M,\lambda)$ as the associated labelled Motzkin path of length $n$.  Then $\pi_{n}$  avoids 321 if and only if $\lambda\,=\,\nu$ (where $\nu$ is the { unitary labelling }) and all horizontal steps in $M$ are at height 0.\\}\\

\textbf{Theorem 8.2.} \textit{Let $\sigma_{n}\in I(321)_{n}$ with $(M,\nu)$ as the associated labelled Motzkin path of length $n$ (whith $\nu$ as before). Then $\sigma_{n}$ is simple if and only if both the following properties hold:\\
i) $(M,\nu)$ is an irreducible  Dyck path;\\
ii) Let $\left\{U_{1},\ldots,U_{m}\right\}$ and  $\left\{D_{1},\ldots,D_{m}\right\}$ be the sequences of the up and of the down steps in $(M,\nu)$. If two up steps $U_{i}$ and $U_{i+1}$ are consecutive up steps in $(M,\nu)$,  then the corresponding $D_{i}$ and $D_{i+1}$ are never consecutive down steps in $(M,\nu)$. }\\

\begin{proof}  Let  $\sigma_{n}\in I(321)_{n}$ be simple, so with no fixed points: then by construction $(M,\nu)$  has no horizontal steps, therefore being a Dyck path, irreducible because a simple involution is connected.\\
Moreover $(M,\nu)$ is such that, always by construction, the maxima of $\sigma$ correspond to the up steps, the minima to the down steps.  Hence, by Proposition 6.2, if $U_{i},\, U_{i+1}$ are consecutive in $(M,\nu)$, $D_{i},\, D_{i+1}$ cannot be consecutive.\\
Conversely, if $(M,\nu)$ satisfies i) and ii), the involution $\sigma\in I(321)_{n}$ is connected, so not of type 12. It follows from Theorem 6.3 that $\sigma_{n}$ is neither of type 21, nor an expansion of a simple involution. Therefore $\sigma$ is simple, as requested.
\end{proof}

In such a way we have shown, as claimed, a second procedure to determine if an involution of $I(321)$ is simple.\\
{Note that for irreducible involutions avoiding $(321)$, the Motzkin path unitary labelled reduces to a Dyck path, leading to the same bijection exposed in \cite{M}. }\\

\noindent Fig. \ref{Figure 5} shows the Dyck path corresponding to an involution of type 12, i.e. sum decomposable : the claim \textit{i)} is not satisfied and {the irreducible components can be recognized in the picture}.\\ Fig. \ref{Figure 6} shows the  Dyck path of an inflation of a simple involution of length greater than 2: { in this case the assumption i) holds, since $\sigma $ is not of type 12, but} \textit{ii)} is not true because the up steps 2,3 correspond to the down steps 7,8  and the up steps 5,6 to the down steps 9,10.

\begin{figure}[h]

\begin{center}
 \setlength{\unitlength}{2mm}
 \begin{picture}(30,6)

 \put(0,0){\begin{picture} (0,0)
 \put(0,0){\circle*{0.3}} \put(2,2){\circle*{0.3}} \put(4,4){\circle*{0.3}} \put(6,2){\circle*{0.3}} \put(8,4){\circle*{0.3}} \put(10,2){\circle*{0.3}} \put(12,0){\circle*{0.3}} \put(14,2){\circle*{0.3}} \put(16,4){\circle*{0.3}} \put(18,2){\circle*{0.3}} \put(20,4){\circle*{0.3}} \put(22,2){\circle*{0.3}}
 \put(24,0){\circle*{0.3}}  \put(26,2){\circle*{0.3}}  \put(28,4){\circle*{0.3}}  \put(30,2){\circle*{0.3}}  \put(32,4){\circle*{0.3}} \put(34,2){\circle*{0.3}} \put(36,0){\circle*{0.3}}

\put(0,0){\line(1,1){4}} \put(4,4){\line(1,-1){2}} \put(6,2){\line(1,1){2}} \put(8,4){\line(1,-1){4}} \put(12,0){\line(1,1){4}} \put(16,4){\line(1,-1){2}} \put(18,2){\line(1,1){2}} \put(20,4){\line(1,-1){4}}
\put(24,0){\line(1,1){4}} \put(28,4){\line(1,-1){2}} \put(30,2){\line(1,1){2}} \put(32,4){\line(1,-1){4}}
 \end{picture}}
 \end{picture}
 \end{center}
 \caption{Example of involution of type 12 (the one of Fig.\ref{Figure 4}) \label{Figure 5}}
 \end{figure}
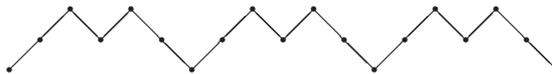

\begin{figure}[h]

\begin{center}
 \setlength{\unitlength}{2mm}
 \begin{picture}(30,6)

 \put(4,0){\begin{picture} (0,0)
 \put(0,0){\circle*{0.3}} \put(2,2){\circle*{0.3}} \put(4,4){\circle*{0.3}} \put(6,6){\circle*{0.3}} \put(8,4){\circle*{0.3}} \put(10,6){\circle*{0.3}} \put(12,8){\circle*{0.3}} \put(14,6){\circle*{0.3}} \put(16,4){\circle*{0.3}} \put(18,2){\circle*{0.3}} \put(20,0){\circle*{0.3}}

\put(0,0){\line(1,1){6}} \put(6,6){\line(1,-1){2}} \put(8,4){\line(1,1){4}} \put(12,8){\line(1,-1){8}}
 \end{picture}}
 \end{picture}
 \end{center}
 \caption{Dyck path corresponding to $\,351624\,[1,\,12,\,1,\,12,\,12,\,12]$ \label{Figure 6}}
 \end{figure}
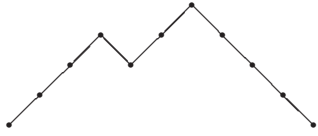

\section{Appendix}

In Section 7, Fig. \ref{Figure 2} and Fig. \ref{Figure 3}, we studied the simple involutions of length 6 and 8, their plots, their sequences and their short Motzkin paths.\\ Now we consider the simple involutions for $n=10$.\\
$
\textbf{35}1\textbf{7}2 \textbf{9}4 (\textbf{10}) 68\, =\, (31)(52)(74)(96)((10)8),$ with $\{1,3,3,3,1\}$ as sequence;\\
$\textbf{36}1\textbf{79}24 (\textbf{10})58\,=\,(31)(62)(74)(95)((10)8)$, corresponding to $\{1,3,1,3,1\}$\\
$\textbf{468}1\textbf{9}2 (\textbf{10})357\,=\,(41)(62)(83)(95)((10)7)$ corresponding to $\{1,3,5,3,1\}$.\\

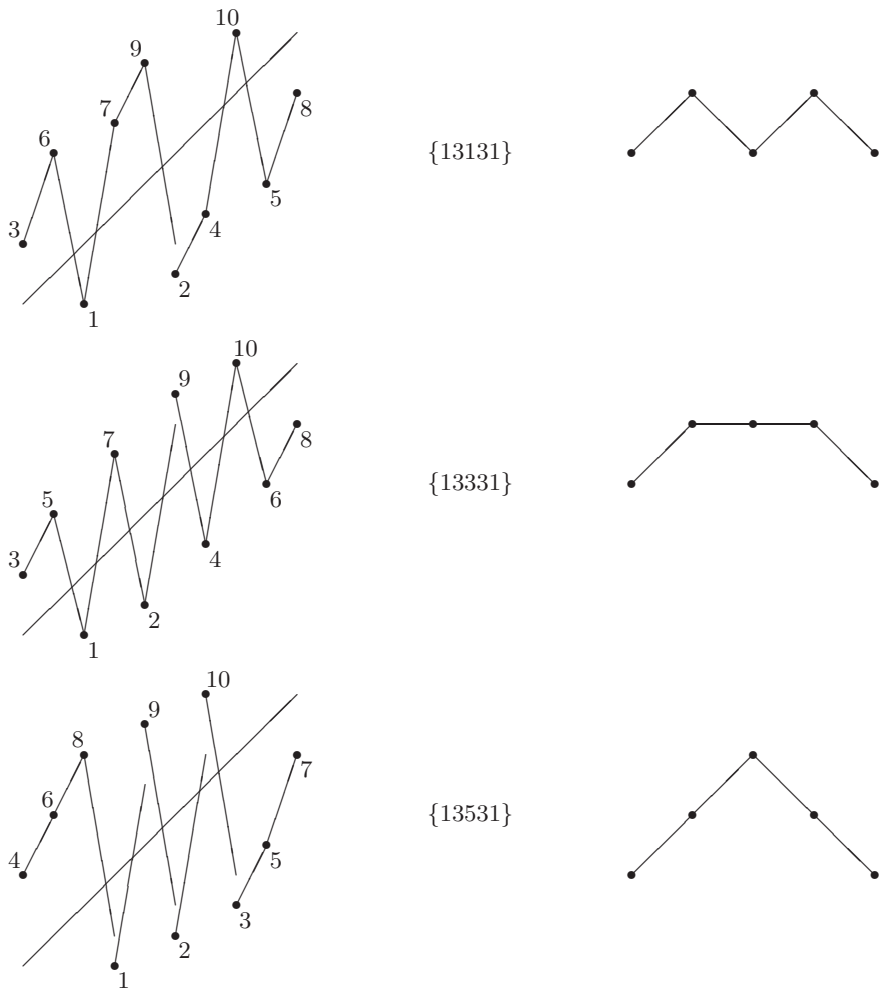
\begin{figure}[h]

\begin{center}
 \setlength{\unitlength}{4mm}
 \begin{picture}(30,10)

\put(0,0){
\begin {picture}(0,0)
\put(0,2){\circle*{0.3}} \put(1,5){\circle*{0.3}} \put(2,0){\circle*{0.3}}\put(3,6){\circle*{0.3}} \put(5,1){\circle*{0.3}} \put(4,8){\circle*{0.3}} \put(6,3){\circle*{0.3}} \put(7,9){\circle*{0.3}} \put(8,4){\circle*{0.3}} \put(9,7){\circle*{0.3}}

\put(0,2){\line(1,3){1}} \put(1,5){\line(1,-5){1}} \put(2,0){\line(1,6){1}} \put(3,6){\line(1,2){1}} \put(4,8){\line(1,-6){1}}
\put(5,1){\line(1,2){1}} \put(6,3){\line(1,6){1}} \put(7,9){\line(1,-5){1}} \put(8,4){\line(1,3){1}}
\put(0,0){\line(1,1){9}}

\put(-0.3,2.5){\makebox(0,0){\small $3$}}
\put(0.7,5.5){\makebox(0,0){\small $6$}}
\put(2.3,-0.5){\makebox(0,0){\small $1$}}
\put(2.7,6.5){\makebox(0,0){\small $7$}}
\put(5.3,0.5){\makebox(0,0){\small $2$}}
\put(3.7,8.5){\makebox(0,0){\small $9$}}
\put(6.7,9.5){\makebox(0,0){\small $10$}}
\put(8.3,3.5){\makebox(0,0){\small $5$}}
\put(9.3,6.5){\makebox(0,0){\small $8$}}

\put(6.3,2.5){\makebox(0,0){\small $4$}}

\end{picture}}

\put(15,5){\makebox(0,0){\small $\{13131\}$}}

 \put(20,5){
 \begin{picture}(0,0)
 \put(0,0){\circle*{0.3}} \put(2,2){\circle*{0.3}} \put(4,0){\circle*{0.3}} \put(6,2){\circle*{0.3}} \put(8,0){\circle*{0.3}}
\put(0,0){\line(1,1){2}} \put(2,2){\line(1,-1){2}} \put(4,0){\line(1,1){2}}\put(6,2){\line(1,-1){2}}
 \end{picture}}
\end{picture}

\end{center}

\begin{center}
 \setlength{\unitlength}{4mm}
 \begin{picture}(30,10)

\put(0,0){
\begin {picture}(0,0)
\put(0,2){\circle*{0.3}} \put(1,4){\circle*{0.3}} \put(2,0){\circle*{0.3}}\put(3,6){\circle*{0.3}} \put(4,1){\circle*{0.3}} \put(5,8){\circle*{0.3}} \put(6,3){\circle*{0.3}} \put(7,9){\circle*{0.3}} \put(8,5){\circle*{0.3}} \put(9,7){\circle*{0.3}}
\put(0,2){\line(1,2){1}} \put(1,4){\line(1,-4){1}} \put(2,0){\line(1,6){1}} \put(3,6){\line(1,-5){1}} \put(4,1){\line(1,6){1}}
\put(5,8){\line(1,-5){1}} \put(6,3){\line(1,6){1}} \put(7,9){\line(1,-4){1}}\put(8,5){\line(1,2){1}}
\put(0,0){\line(1,1){9}}

\put(-0.3,2.5){\makebox(0,0){\small $3$}}
\put(0.8,4.5){\makebox(0,0){\small $5$}}
\put(2.3,-0.5){\makebox(0,0){\small $1$}}
\put(2.8,6.5){\makebox(0,0){\small $7$}}
\put(4.3,0.5){\makebox(0,0){\small $2$}}
\put(5.3,8.5){\makebox(0,0){\small $9$}}
\put(6.3,2.5){\makebox(0,0){\small $4$}}
\put(8.3,4.5){\makebox(0,0){\small $6$}}
\put(7.3,9.5){\makebox(0,0){\small $10$}}
\put(9.3,6.5){\makebox(0,0){\small $8$}}

\end{picture}}

\put(15,5){\makebox(0,0){\small $\{13331\}$}}

 \put(20,5){
 \begin{picture}(0,0)
  \put(0,0){\circle*{0.3}} \put(2,2){\circle*{0.3}} \put(4,2){\circle*{0.3}} \put(6,2){\circle*{0.3}} \put(8,0){\circle*{0.3}}
\put(0,0){\line(1,1){2}} \put(2,2){\line(1,0){4}} \put(6,2){\line(1,-1){2}}
 \end{picture}}
\end{picture}

\end{center}

\begin{center}
 \setlength{\unitlength}{4mm}
 \begin{picture}(30,10)

\put(0,0){
\begin {picture}(0,0)
\put(0,3){\circle*{0.3}} \put(1,5){\circle*{0.3}} \put(2,7){\circle*{0.3}}
\put(3,0){\circle*{0.3}}\put(4,8){\circle*{0.3}} \put(5,1){\circle*{0.3}} \put(6,9){\circle*{0.3}} \put(7,2){\circle*{0.3}} \put(8,4){\circle*{0.3}} \put(9,7){\circle*{0.3}}

\put(0,3){\line(1,2){1}} \put(1,5){\line(1,2){1}} \put(2,7){\line(1,-6){1}} \put(3,0){\line(1,6){1}} \put(4,8){\line(1,-6){1}}
\put(5,1){\line(1,6){1}} \put(6,9){\line(1,-6){1}} \put(7,2){\line(1,2){1}} \put(8,4){\line(1,3){1}}
\put(0,0){\line(1,1){9}}

\put(-0.3,3.5){\makebox(0,0){\small $4$}}
\put(0.8,5.5){\makebox(0,0){\small $6$}}
\put(1.8,7.5){\makebox(0,0){\small $8$}}
\put(4.3,8.5){\makebox(0,0){\small $9$}}
\put(6.4,9.5){\makebox(0,0){\small $10$}}
\put(3.3,-0.5){\makebox(0,0){\small $1$}}
\put(5.3,0.5){\makebox(0,0){\small $2$}}
\put(7.3,1.5){\makebox(0,0){\small $3$}}
\put(9.3,6.5){\makebox(0,0){\small $7$}}
\put(8.3,3.5){\makebox(0,0){\small $5$}}

\end{picture}}

\put(15,5){\makebox(0,0){\small $\{13531\}$}}

 \put(20,3){
 \begin{picture}(0,0)
  \put(0,0){\circle*{0.3}} \put(2,2){\circle*{0.3}} \put(4,4){\circle*{0.3}} \put(6,2){\circle*{0.3}} \put(8,0){\circle*{0.3}}
\put(0,0){\line(1,1){4}} \put(4,4){\line(1,-1){4}}
 \end{picture}}
\end{picture}
\end{center}

\end{figure}

{Finally we list the simple involutions avoiding $(321)$ for $n=12$ and $n=14$ together with their sequences and their Motzkin path of length 5 and 6.
\begin{figure}[h]

\begin{center}
 \setlength{\unitlength}{4mm}
 \begin{picture}(30,2)

 \put(0,0){\makebox (0,0) {\small \textbf{3 5}1 \textbf{7}2 \textbf{9}4 (\textbf{11})6 (\textbf{12})8(10),}}
 \put(12,0){\makebox (0,0) {\small \{133331\}}}
 \put(20,0){\begin{picture} (0,0)
 \put(0,0){\circle*{0.3}} \put(1,1){\circle*{0.3}} \put(2,1){\circle*{0.3}} \put(3,1){\circle*{0.3}} \put(4,1){\circle*{0.3}} \put(5,0){\circle*{0.3}}

\put(0,0){\line(1,1){1}} \put(1,1){\line(1,0){3}} \put(4,1){\line(1,-1){1}}
 \end{picture}}
 \end{picture}
 \setlength{\unitlength}{4mm}
 \begin{picture}(30,2)

 \put(0,0){\makebox (0,0) {\small \textbf{3 5}1 \textbf{8}2 \textbf{9 11}46 (\textbf{12})7(10),}}
 \put(12,0){\makebox (0,0) {\small \{133131\}}}
 \put(20,0){\begin{picture} (0,0)
 \put(0,0){\circle*{0.3}} \put(1,1){\circle*{0.3}} \put(2,1){\circle*{0.3}} \put(3,0){\circle*{0.3}} \put(4,1){\circle*{0.3}} \put(5,0){\circle*{0.3}}

\put(0,0){\line(1,1){1}} \put(1,1){\line(1,0){1}} \put(2,1){\line(1,-1){1}} \put(3,0){\line(1,1){1}} \put(4,1){\line(1,-1){1}}
 \end{picture}}
 \end{picture}

 \setlength{\unitlength}{4mm}
 \begin{picture}(30,2)

 \put(0,0){\makebox (0,0) {\small \textbf{3 6}1 \textbf{7 9}24 (\textbf{11})5 (\textbf{12})8(10),}}
 \put(12,0){\makebox (0,0) {\small \{131331\}}}
 \put(20,0){\begin{picture} (0,0)
 \put(0,0){\circle*{0.3}} \put(1,1){\circle*{0.3}} \put(2,0){\circle*{0.3}} \put(3,1){\circle*{0.3}} \put(4,1){\circle*{0.3}} \put(5,0){\circle*{0.3}}

\put(0,0){\line(1,1){1}} \put(1,1){\line(1,-1){1}} \put(2,0){\line(1,1){1}} \put(3,1){\line(1,0){1}} \put(4,1){\line(1,-1){1}}
 \end{picture}}
 \end{picture}

  \setlength{\unitlength}{4mm}
 \begin{picture}(30,3)

 \put(0,0){\makebox (0,0) {\small \textbf{3 6}1 \textbf{8 (10)}2 (\textbf{11})4 (\textbf{12})579,}}
 \put(12,0){\makebox (0,0) {\small \{133531\}}}
 \put(20,0){\begin{picture} (0,0)
 \put(0,0){\circle*{0.3}} \put(1,1){\circle*{0.3}} \put(2,1){\circle*{0.3}} \put(3,2){\circle*{0.3}} \put(4,1){\circle*{0.3}} \put(5,0){\circle*{0.3}}

\put(0,0){\line(1,1){1}} \put(1,1){\line(1,0){1}} \put(2,1){\line(1,1){1}} \put(3,2){\line(1,-1){1}} \put(4,1){\line(1,-1){1}}
 \end{picture}}
 \end{picture}

 \setlength{\unitlength}{4mm}
 \begin{picture}(30,3)

 \put(0,0){\makebox (0,0) {\small \textbf{4 6 8}1 \textbf{9}2 (\textbf{11})35 (\textbf{12})7(10),}}
 \put(12,0){\makebox (0,0) {\small \{135331\}}}
 \put(20,0){\begin{picture} (0,0)
 \put(0,0){\circle*{0.3}} \put(1,1){\circle*{0.3}} \put(2,2){\circle*{0.3}} \put(3,1){\circle*{0.3}} \put(4,1){\circle*{0.3}} \put(5,0){\circle*{0.3}}

\put(0,0){\line(1,1){2}} \put(2,2){\line(1,-1){1}} \put(3,1){\line(1,0){1}} \put(4,1){\line(1,-1){1}}
 \end{picture}}
 \end{picture}

  \setlength{\unitlength}{4mm}
 \begin{picture}(30,3)

 \put(0,0){\makebox (0,0) {\small \textbf{4 6 8}1 (\textbf{10})2 (\textbf{11})3 (\textbf{12})579.}}
 \put(12,0){\makebox (0,0) {\small \{135331\}}}
 \put(20,0){\begin{picture} (0,0)
 \put(0,0){\circle*{0.3}} \put(1,1){\circle*{0.3}} \put(2,2){\circle*{0.3}} \put(3,2){\circle*{0.3}} \put(4,1){\circle*{0.3}} \put(5,0){\circle*{0.3}}

\put(0,0){\line(1,1){2}} \put(2,2){\line(1,0){1}} \put(3,2){\line(1,-1){1}} \put(4,1){\line(1,-1){1}}
 \end{picture}}
 \end{picture}
 \caption{Simple involutions in $I(321)_{12}$}

 \end{center}
 \end{figure}

\newpage

\begin{figure}[h]

\begin{center}
 \setlength{\unitlength}{4mm}
 \begin{picture}(30,2)

 \put(0,0){\makebox (0,0) {\small \textbf{35}1 \textbf{7}2 \textbf{9}4 \textbf{(11)}6 \textbf{(13)}8 \textbf{14}(10)(12),}}
 \put(12,0){\makebox (0,0) {\small \{1333331\}}}
 \put(20,0){\begin{picture} (0,0)
 \put(0,0){\circle*{0.3}} \put(1,1){\circle*{0.3}} \put(2,1){\circle*{0.3}} \put(3,1){\circle*{0.3}} \put(4,1){\circle*{0.3}} \put(5,1){\circle*{0.3}} \put(6,0){\circle*{0.3}}

\put(0,0){\line(1,1){1}} \put(1,1){\line(1,0){4}} \put(5,1){\line(1,-1){1}}
 \end{picture}}
 \end{picture}

  \setlength{\unitlength}{4mm}
 \begin{picture}(30,2)

 \put(0,0){\makebox (0,0) {\small \textbf{35}1 \textbf{7}2 \textbf{(10)}4 \textbf{(11)(13)}68 \textbf{(14)}9(12),}}
 \put(12,0){\makebox (0,0) {\small \{1333131\}}}
 \put(20,0){\begin{picture} (0,0)
 \put(0,0){\circle*{0.3}} \put(1,1){\circle*{0.3}} \put(2,1){\circle*{0.3}} \put(3,1){\circle*{0.3}} \put(4,0){\circle*{0.3}} \put(5,1){\circle*{0.3}} \put(6,0){\circle*{0.3}}

\put(0,0){\line(1,1){1}} \put(1,1){\line(1,0){2}} \put(3,1){\line(1,-1){1}} \put(4,0){\line(1,1){1}} \put(5,1){\line(1,-1){1}}
 \end{picture}}
 \end{picture}

  \setlength{\unitlength}{4mm}
 \begin{picture}(30,2)

 \put(0,0){\makebox (0,0) {\small \textbf{35}1 \textbf{8}2 \textbf{9 (11)}46 \textbf{(13)}7 \textbf{(14)}(10)(12),}}
 \put(12,0){\makebox (0,0) {\small \{1331331\}}}
 \put(20,0){\begin{picture} (0,0)
 \put(0,0){\circle*{0.3}} \put(1,1){\circle*{0.3}} \put(2,1){\circle*{0.3}} \put(3,0){\circle*{0.3}} \put(4,1){\circle*{0.3}} \put(5,1){\circle*{0.3}} \put(6,0){\circle*{0.3}}

\put(0,0){\line(1,1){1}} \put(1,1){\line(1,0){1}} \put(2,1){\line(1,-1){1}} \put(3,0){\line(1,1){1}} \put(4,1){\line(1,0){1}} \put(5,1){\line(1,-1){1}}
 \end{picture}}
 \end{picture}

   \setlength{\unitlength}{4mm}
 \begin{picture}(30,2)

 \put(0,0){\makebox (0,0) {\small \textbf{36}1 \textbf{79}2 4\textbf{(11)}5 \textbf{(13)}8 \textbf{(14)}(10)(12),}}
 \put(12,0){\makebox (0,0) {\small \{1313331\}}}
 \put(20,0){\begin{picture} (0,0)
 \put(0,0){\circle*{0.3}} \put(1,1){\circle*{0.3}} \put(2,0){\circle*{0.3}} \put(3,1){\circle*{0.3}} \put(4,1){\circle*{0.3}} \put(5,1){\circle*{0.3}} \put(6,0){\circle*{0.3}}

\put(0,0){\line(1,1){1}} \put(1,1){\line(1,-1){1}} \put(2,0){\line(1,1){1}} \put(3,1){\line(1,0){1}} \put(4,1){\line(1,0){1}} \put(5,1){\line(1,-1){1}}
 \end{picture}}
 \end{picture}

 \setlength{\unitlength}{4mm}
 \begin{picture}(30,2)
 \put(0,0){\makebox (0,0) {\small \textbf{36}1 \textbf{7(10)}2 4\textbf{(11)(13)}5 8 \textbf{(14)}9(12),}}
 \put(12,0){\makebox (0,0) {\small \{1313131\}}}
 \put(20,0){\begin{picture} (0,0)
 \put(0,0){\circle*{0.3}} \put(1,1){\circle*{0.3}} \put(2,0){\circle*{0.3}} \put(3,1){\circle*{0.3}} \put(4,0){\circle*{0.3}} \put(5,1){\circle*{0.3}} \put(6,0){\circle*{0.3}}

\put(0,0){\line(1,1){1}} \put(1,1){\line(1,-1){1}} \put(2,0){\line(1,1){1}} \put(3,1){\line(1,-1){1}} \put(4,0){\line(1,1){1}} \put(5,1){\line(1,-1){1}}

 \end{picture}}
 \end{picture}

  \setlength{\unitlength}{4mm}
 \begin{picture}(30,3)

 \put(0,0){\makebox (0,0) {\small \textbf{35}1 \textbf{8}2 \textbf{(10)(12)}4 \textbf{(13)}6 \textbf{(14)}79(11),}}
 \put(12,0){\makebox (0,0) {\small \{1333531\}}}
 \put(20,0){\begin{picture} (0,0)
 \put(0,0){\circle*{0.3}} \put(1,1){\circle*{0.3}} \put(2,1){\circle*{0.3}} \put(3,1){\circle*{0.3}} \put(4,2){\circle*{0.3}} \put(5,1){\circle*{0.3}} \put(6,0){\circle*{0.3}}

\put(0,0){\line(1,1){1}} \put(1,1){\line(1,0){1}} \put(2,1){\line(1,0){1}} \put(3,1){\line(1,1){1}} \put(4,2){\line(1,-1){1}} \put(5,1){\line(1,-1){1}}
 \end{picture}}
 \end{picture}

 \setlength{\unitlength}{4mm}
 \begin{picture}(30,3)

 \put(0,0){\makebox (0,0) {\small \textbf{36}1 \textbf{8(10)}2 \textbf{(12)}4 \textbf{(13)}5 \textbf{(14)}79(11),}}
 \put(12,0){\makebox (0,0) {\small \{1335531\}}}
 \put(20,0){\begin{picture} (0,0)
 \put(0,0){\circle*{0.3}} \put(1,1){\circle*{0.3}} \put(2,1){\circle*{0.3}} \put(3,2){\circle*{0.3}} \put(4,2){\circle*{0.3}} \put(5,1){\circle*{0.3}} \put(6,0){\circle*{0.3}}

\put(0,0){\line(1,1){1}} \put(1,1){\line(1,0){1}} \put(2,1){\line(1,1){1}} \put(3,2){\line(1,0){1}} \put(4,2){\line(1,-1){1}} \put(5,1){\line(1,-1){1}}
 \end{picture}}
 \end{picture}

 \setlength{\unitlength}{4mm}
 \begin{picture}(30,3)

 \put(0,0){\makebox (0,0) {\small \textbf{36}1 \textbf{8(10)}2 \textbf{(11)}4 \textbf{(13)}57 \textbf{(14)}9(12),}}
 \put(12,0){\makebox (0,0) {\small \{1335331\}}}
 \put(20,0){\begin{picture} (0,0)
 \put(0,0){\circle*{0.3}} \put(1,1){\circle*{0.3}} \put(2,1){\circle*{0.3}} \put(3,2){\circle*{0.3}} \put(4,1){\circle*{0.3}} \put(5,1){\circle*{0.3}} \put(6,0){\circle*{0.3}}

\put(0,0){\line(1,1){1}} \put(1,1){\line(1,0){1}} \put(2,1){\line(1,1){1}} \put(3,2){\line(1,-1){1}} \put(4,1){\line(1,0){1}} \put(5,1){\line(1,-1){1}}
 \end{picture}}
 \end{picture}

\setlength{\unitlength}{4mm}
 \begin{picture}(30,3)

 \put(0,0){\makebox (0,0) {\small \textbf{468}1 \textbf{9}2 \textbf{(11)}35 \textbf{(13)}7 \textbf{(14)}(10)(12),}}
 \put(12,0){\makebox (0,0) {\small \{1353331\}}}
 \put(20,0){\begin{picture} (0,0)
 \put(0,0){\circle*{0.3}} \put(1,1){\circle*{0.3}} \put(2,2){\circle*{0.3}} \put(3,1){\circle*{0.3}} \put(4,1){\circle*{0.3}} \put(5,1){\circle*{0.3}} \put(6,0){\circle*{0.3}}

\put(0,0){\line(1,1){1}} \put(1,1){\line(1,1){1}} \put(2,2){\line(1,-1){1}} \put(3,1){\line(1,0){1}} \put(4,1){\line(1,0){1}} \put(5,1){\line(1,-1){1}}
 \end{picture}}
 \end{picture}

 \setlength{\unitlength}{4mm}
 \begin{picture}(30,3)

 \put(0,0){\makebox (0,0) {\small \textbf{37}1 \textbf{8(10)(12)}24 \textbf{(13)}5 \textbf{(14)}69(11),}}
 \put(12,0){\makebox (0,0) {\small \{1313531\}}}
 \put(20,0){\begin{picture} (0,0)
 \put(0,0){\circle*{0.3}} \put(1,1){\circle*{0.3}} \put(2,0){\circle*{0.3}} \put(3,1){\circle*{0.3}} \put(4,2){\circle*{0.3}} \put(5,1){\circle*{0.3}} \put(6,0){\circle*{0.3}}

\put(0,0){\line(1,1){1}} \put(1,1){\line(1,-1){1}} \put(2,0){\line(1,1){1}} \put(3,1){\line(1,1){1}} \put(4,2){\line(1,-1){1}} \put(5,1){\line(1,-1){1}}
 \end{picture}}
 \end{picture}

  \setlength{\unitlength}{4mm}
 \begin{picture}(30,3)

 \put(0,0){\makebox (0,0) {\small \textbf{469}1 \textbf{(10)}2 \textbf{(11)(13)}357 \textbf{(14)}8(12),}}
 \put(12,0){\makebox (0,0) {\small \{1353131\}}}
 \put(20,0){\begin{picture} (0,0)
 \put(0,0){\circle*{0.3}} \put(1,1){\circle*{0.3}} \put(2,2){\circle*{0.3}} \put(3,1){\circle*{0.3}} \put(4,0){\circle*{0.3}} \put(5,1){\circle*{0.3}} \put(6,0){\circle*{0.3}}

\put(0,0){\line(1,1){1}} \put(1,1){\line(1,1){1}} \put(2,2){\line(1,-1){1}} \put(3,1){\line(1,-1){1}} \put(4,0){\line(1,1){1}} \put(5,1){\line(1,-1){1}}
 \end{picture}}
 \end{picture}

 \setlength{\unitlength}{4mm}
 \begin{picture}(30,3)

 \put(0,0){\makebox (0,0) {\small \textbf{468}1 \textbf{(10)}2 \textbf{(11)}3 \textbf{(13)}57 \textbf{(14)}9(12),}}
 \put(12,0){\makebox (0,0) {\small \{1355331\}}}
 \put(20,0){\begin{picture} (0,0)
 \put(0,0){\circle*{0.3}} \put(1,1){\circle*{0.3}} \put(2,2){\circle*{0.3}} \put(3,2){\circle*{0.3}} \put(4,1){\circle*{0.3}} \put(5,1){\circle*{0.3}} \put(6,0){\circle*{0.3}}

\put(0,0){\line(1,1){1}} \put(1,1){\line(1,1){1}} \put(2,2){\line(1,0){1}} \put(3,2){\line(1,-1){1}} \put(4,1){\line(1,0){1}} \put(5,1){\line(1,-1){1}}
 \end{picture}}
 \end{picture}

  \setlength{\unitlength}{4mm}
 \begin{picture}(30,3)

 \put(0,0){\makebox (0,0) {\small \textbf{468}1 \textbf{(10)}2 \textbf{(12)}3 \textbf{(13)}5 \textbf{(14)}79(11),}}
 \put(12,0){\makebox (0,0) {\small \{1355531\}}}
 \put(20,0){\begin{picture} (0,0)
 \put(0,0){\circle*{0.3}} \put(1,1){\circle*{0.3}} \put(2,2){\circle*{0.3}} \put(3,2){\circle*{0.3}} \put(4,2){\circle*{0.3}} \put(5,1){\circle*{0.3}} \put(6,0){\circle*{0.3}}

\put(0,0){\line(1,1){1}} \put(1,1){\line(1,1){1}} \put(2,2){\line(1,0){1}} \put(3,2){\line(1,0){1}} \put(4,2){\line(1,-1){1}} \put(5,1){\line(1,-1){1}}
 \end{picture}}
 \end{picture}

 \setlength{\unitlength}{4mm}
 \begin{picture}(30,3)

 \put(0,0){\makebox (0,0) {\small \textbf{479}1 \textbf{(10)(12)}2 \textbf{(13)}35 \textbf{(14)}68(11),}}
 \put(12,0){\makebox (0,0) {\small \{1353531\}}}
 \put(20,0){\begin{picture} (0,0)
 \put(0,0){\circle*{0.3}} \put(1,1){\circle*{0.3}} \put(2,2){\circle*{0.3}} \put(3,1){\circle*{0.3}} \put(4,2){\circle*{0.3}} \put(5,1){\circle*{0.3}} \put(6,0){\circle*{0.3}}

\put(0,0){\line(1,1){1}} \put(1,1){\line(1,1){1}} \put(2,2){\line(1,-1){1}} \put(3,1){\line(1,1){1}} \put(4,2){\line(1,-1){1}} \put(5,1){\line(1,-1){1}}
 \end{picture}}
 \end{picture}

  \setlength{\unitlength}{4mm}
 \begin{picture}(30,4)

 \put(0,0){\makebox (0,0) {\small \textbf{579(11)}1 \textbf{(12)}2 \textbf{(13)}3 \textbf{(14)}468(10),}}
 \put(12,0){\makebox (0,0) {\small \{1357531\}}}
 \put(20,0){\begin{picture} (0,0)
 \put(0,0){\circle*{0.3}} \put(1,1){\circle*{0.3}} \put(2,2){\circle*{0.3}} \put(3,3){\circle*{0.3}} \put(4,2){\circle*{0.3}} \put(5,1){\circle*{0.3}} \put(6,0){\circle*{0.3}}

\put(0,0){\line(1,1){1}} \put(1,1){\line(1,1){1}} \put(2,2){\line(1,1){1}} \put(3,3){\line(1,-1){1}} \put(4,2){\line(1,-1){1}} \put(5,1){\line(1,-1){1}}
 \end{picture}}
 \end{picture}

\caption{Simple involutions in $I(321)_{14}$}

 \end{center} \end{figure}


\begin{thebibliography}{99}

\bibitem{a}
M.H. Albert, \emph{The fine structure of 321 avoiding permutations}, Technical Report OUCS-2002-11
\bibitem{AA}
M.H. Albert, M.D. Atkinson, \emph{Simple permutations and pattern restricted permutations}, Discrete Mathematics 300 (2005) 1 - 15
\bibitem{AAA}
M. Barnabei et al., \emph{Restricted involutions and Motzkin paths}, Adv. in Appl. Math. (2010), doi:10.1016/j.aam.2010.05.002
\bibitem{B}
 R. Brignall, S. Huczynska, V. Vatter, \emph{Simple permutations and algebraic generating functions}, J. of Combinatorial Theory, Series A115(2008) 423-441
 \bibitem{MR}
 R.H. M\"{o}hring, F.J. Radermacher, \emph{Substitution decomposition for discrete structures and connections with combinatorial optimization}, North-Holland Math.Stud.vol.95 (1984), 257 - 355
 \bibitem{M}
E. Munarini, C. Perelli Cippo, \emph{Statistics and codes on linear permutations}, Permutation Patterns 2009.
\bibitem{ss}
R. Simion, F.W. Schmidt, \emph{Restricted permutations}, Europ.T.Combin. 6 (1985), 383 - 406
\bibitem{S}
  N. J. A. Sloane,
  \emph{On-Line Encyclopedia of Integer Sequences},
  published electronically at {\tt http://www.research.att.com/\~{}njas/sequences/}.



\end{thebibliography}
\end{document}